\documentclass{article}%
\usepackage{amsmath}
\usepackage{graphicx}
\usepackage{amsfonts}
\usepackage{amssymb}%
\newtheorem{theorem}{Theorem}

\newtheorem{lemma}[theorem]{Lemma}

\newtheorem{proposition}[theorem]{Proposition}

\newenvironment{proof}[1][Proof]{\noindent\textbf{#1.} }{\ \rule{0.5em}{0.5em}}
\begin{document}

\title{Classification of Steadily Rotating Spiral Waves for the Kinematic
Model\thanks{This work was supported in part by the National Science Council
of Taiwan and the Natural Sciences and Engineering Research Council of
Canada.}}
\author{Chu-Pin Lo\thanks{Department of Applied Mathematics, Providence University,
200 Chung-chi Rd., Shalu, Taichung County, Taiwan 433, \texttt{cplo@pu.edu.tw}%
.}
\and Nedialko~S.~Nedialkov\thanks{Department of Computing and Software, McMaster
University, Hamilton, Ontario, Canada, L8S 4L7, \texttt{nedialk@mcmaster.ca}.}
\and Juan-Ming Yuan\thanks{Department of Mathematics and the Texas Institute for
Computational and Applied Mathematics,The University of Texas, Austin, TX
78712, U.S.A., \texttt{jmyuan@math.utexas.edu}.}}
\maketitle

\begin{abstract}
Spiral waves arise in many biological, chemical, and physiological systems.
The kinematical model can be used to describe the motion of the spiral arms
approximated as curves in the plane. For this model, there appeared some
results in the literature. However, these results all are based upon some
simplification on the model or prior phenomenological assumptions on the
solutions. In this paper, we use really full kinematic model to classify a
generic kind of steadily rotating spiral waves, i.e., with positive (or
negative) curvature. In fact, using our results (Theorem \ref{Theorem6}), we
can answer the following questions: Is there any steadily rotating spiral wave
for a given weakly excitable medium? If yes, what kind of information we can
know about these spiral waves? e.g., the tip's curvature, the tip's tangential
velocity, and the rotating frequency. Comparing our results with previous ones
in the literature, there are some differences between them. There are only
solutions with monotonous curvatures via simplified model but full model
admits solutions with any given oscillating number of the curvatures.

\end{abstract}

\pagestyle{myheadings} \thispagestyle{plain}
\markboth{C.~-P.~LO, N.~S.~NEDIALKOV, AND J.-M.~YUAN}{CLASSIFICATION OF
STEADILY ROTATING SPIRAL WAVES}

\textbf{Key Words: }kinematic model, spiral waves, excitable media

\textbf{AMS subject classifications}. 37N25, 47N60, 93A30, 47N20, 65L05, 65G20

\section{Introduction}

Rotating spiral waves appear in many biological, physiological, and chemical
systems. For example, spiral waves arise in cardiac arrytmias
\cite{GaQ,KeS,KeT2,MRA,Rot}, egg fertilization \cite{Dup}, the
Belousov-Zhabotinsky (BZ) chemical reaction \cite{BBS,JSW,
KeT1,KCZ,MPH1,MPH2,ZKR}, catalysis \cite{CEI,NVR}, and aggregation of slime
mold amoeba \cite{AlM,Dur,SHM,Wan}. These waves have been studied extensively
from experimental \cite{BBS,JSW,KCZ,MPH1,MPH2,NMT,SkS,Win1,ZKR}, numerical
\cite{JSW,JaW,Kar,Lug,MiZ,Win2,ZMM}, and analytical
\cite{Bra,Dup,GLT,Hag,IIY,Kee1,KLR,SSW1,SSW2,Sch,TyK,YaN} aspects.

Often, such systems exhibit the so-called ``excitable''
property~\cite{Mik,Zyk}. A spiral wave represents an excited moving spatial
region with two spiral-like, thin phase-transition layers (or interfaces):
wave front and wave back. These layers separate two different phase regions,
excited and refractory (or recovery), and intersect at a tip (see Figure 12.5
in \cite{KeS} or Figures 13, 17, and 22 in \cite{Mer}).

Usually, the appearance of spiral waves is not desirable. For instance, heart
attacks are caused by an abrupt shift from rhythmic pumping to spasmodic
convulsion of the heart. Normally, with each heartbeat, an electric wavefront
propagates across the interconnected muscle fibers, causing them to contract.
However, because of some abnormality in the tissue, this wave can become stuck
and start rotating as a spiral wave \cite{Pet,PeE}. Furthermore, if spiral
waves occur in the cortex, they may lead to epileptic seizures. On the retina
or the visual cortex, they may cause hallucinations \cite{BCG1,BCG2}.

Spiral waves can be created by breaking propagating waves \cite{JSW,PMH,SkS}.
In \cite{PMH}, a circular wave is broken mechanically by ejecting a gentle
blast of air into a small section of the wave. The process by which a broken
wave evolves toward a rotating spiral wave has been addressed in
\cite{DMZ,Fif1,Fif2,MeP,MDZ,Zyk}. As pointed out in \cite{Fif2}, along the
interface separating the excited and the refractory regions, there exists a
point where the normal velocity of the interface changes sign. This induces a
twisting action on the interface motion, and as a result, the interface starts
wrapping around some center of rotation, to form a spiral structure.
Eventually, in an appropriate parameter range, a state of steady rotation is
achieved, during which the spiral tip traces a circular trajectory at constant
angular velocity, and the curvatures of the spiral-shaped interfaces keep
constant in time. Such steadily rotating spiral waves have been observed in
many experiments \cite{MPH3,MPH2,Win1}.

Another way to initiate spiral waves is by applying a spatially graded
perturbation to a medium in the recovery phase after excitation \cite{Win4,
Win3}. Besides steadily rotating spiral waves, non-steady forms of rotation,
known as ``meander'' \cite{Win3}, are also confirmed by experiments
\cite{JSW,PMH,SkS} and numerical simulations \cite{Bar,JSW,JaW,Lug,Zyk}. When
the controlled parameters are changed within some suitable range
\cite{JSW,PMH,SkS}, or due to the anisotropic influences of the media
\cite{Rot}, or the effect of electric field \cite{MGP,MPG}, a transition from
the steadily rotating spiral waves to meandering ones may occur. To study the
above mentioned wave patterns, many models have been proposed, which we can
categorize here as ``{phase field}'' type and ``{sharp interface}'' type
models\footnote{We borrow these two terms from material science.}.

For the phase field type of models, well-known examples are the
reaction-diffusion equations (for example, the FitzHugh-Nagumo model for the
action potential wave in neurons \cite{KeS}), the Oregonator model for the BZ
reaction \cite{FiN,TyF}, the bidomain model~\cite{KeS}, and the chemotactic
models \cite{KelS,Wan}. The bidomain model is a mixture of elliptic and
reaction-diffusion type of equations. It is used to describe the electrical
activities of cardiac tissue. The chemotactic models describe the oriented
movement of cells in response to the concentration gradient of chemical
substances in their environment, such as the aggregation phenomena of slime
mold amoeba.

Using the reaction-diffusion equation, the existence of steadily rotating
spiral waves near a homogeneous steady state has been proved formally by Hagan
\cite{Hag} and later rigorously by Scheel \cite{Sch}. Other more complicated
wave patterns, e.g., meandering and hypermeander were also explored by
\cite{AMN,Bar,SSW1,SSW2,Win2}.

For the sharp interface type of models, some typical examples are the free
boundary models \cite{Fif3,Kee2,KeK,KeL,KLR,PeS1,PeS2,TyK} and the kinematic
models \cite{EBH,EBH2,Mik,MDZ,MiZ,Zyk}. Other models are given in
\cite{TyK,YaN}.

If the two phase transition layers, as mentioned above, are very thin, we can
view them as two plane curves meeting at the tip. Then, the free boundary
model consists of an interface equation, a simplified phase field equation,
and a dispersion relation with some boundary conditions. The interface
equation (13) or (38) in \cite{TyK} is derived from the so-called
eikonal-curvature equation (11) there and governs the motion of a wave front;
the simplified phase field equation (37) in \cite{TyK} is derived from the
original phase field model and affects the velocity appearing in the interface
equation; and the dispersion relation (see section 3.3 in \cite{TyK}) reflects
the interaction between spiral waves (refractory tail effect).

In weakly excitable media, spiral waves rotate around a large circle, the
excited region becomes rather narrow, and the spiral waves are sparse. In this
case, we need to consider the wave front only. Thus, the spiral wave can be
viewed as a single curve, which performs a motion in the plane regardless of
the interaction between the spiral waves. That is, the dispersion relation can
be ignored.

The kinematic model is then formulated in terms of a motion of a single curve
with free end; see section 2 for the details. The validity of this kinematic
model has been verified, since various wave patterns appearing in the
experiments can also be produced with this approach \cite{MDZ}.

Finally, we should point out that many main governing equations of the sharp
interface models can be derived from the phase field models by taking various
types of singular limits, e.g., the small parameters may represent the width
of the thin transition layers. Therefore, it is natural to expect that there
should be a close relation between the solutions of these two type of models.
Indeed, it has been found that there is a good agreement for temporal periods
of rotating spiral waves, when comparing the reaction-diffusion approach with
the kinematic one \cite{MiZ}.

In this paper, we use the really \textit{\textquotedblleft
full\textquotedblright\ }kinematic model to prove the existence of steadily
rotating spiral waves. In fact, using our results (Theorem \ref{Theorem6}), we
can answer the following questions: Is there any steadily rotating spiral wave
for a given weakly excitable medium (i.e., $V_{0},D$ are given below). If yes,
what kind of informations we can know about these spiral waves? e.g., the
tip's curvature ($\kappa_{0}$ below), the tip's tangential velocity ($G$
below), and the rotating frequency ($\omega$ below). In the literature, there
appeared some results about kinematic model. However, these results all are
based upon some simplifications on the model or prior assumptions on the
solutions. For examples, \cite{EBH}, \cite{EBH2}, \cite{GLT}, and \cite{MDZ}
considered only simplified model\ (see the Remarks in Section 2); \cite{GLT},
\cite{IIY}, \cite{MDZ}, and \cite{YaN} made a prior assumption on the
solution, i.e., tip's tangential velocity $G$ is $0,$ which was proven to be a
very special case of steadily rotating spiral waves in \cite{EBH}; \cite{IIY}
is valid only for nonrotating ones, i.e., angular frequency $\omega$ is $0$
(see the explanation after (\ref{2.11}) below); \cite{IIY} considered
solutions with positive curvature $\kappa.$ In this paper, we also consider
solutions with unchangeable curvature's sign (i.e., with positive or negative
curvature). However, this condition is generic. Note that the formula (3.11)
in \cite{YaN} is $\omega(0;0,l(0))$ in Theorem \ref{Theorem6}, i.e., a special
case of our results with zero tip's tangential velocity ($G=0$) and monotonous
curvature ($i=0$). Note also that steadily rotating spiral waves with positive
curvatures must have \textit{monotonous} curvatures via simplified model (see
\cite{GLT}); however, using full model in this paper we have confirmed the
existence of \textit{oscillating} ones, i.e., the case \textquotedblleft%
$i\neq0\textquotedblright$ in Theorem \ref{Theorem6}.

In section 2, we recall the kinematic model following Mikhailov \textit{et
al}. \cite{MDZ} and a new version by \cite{EBH}, \cite{EBH2}. In section 3, we
state and then prove the main theorem of this paper. Numerical results are
presented in section 4.

\section{The Kinematic Model}

The kinematic model, formulated in terms of a motion of curves, was first
proposed for describing excitation waves in a cardiac muscle by Wiener and
Rosenblueth~\cite{WiR}. They assumed that a plane curve moves in a normal
direction with a constant velocity. They showed that such a curve, rotating
around an obstacle, forms a spiral (which represents an involute of this
obstacle) and approaches an Archimedian spiral far from it.

Since there occurs a singularity in the curvature of this spiral curve near
the tip \cite{WiR}, and also by the subsequent analysis of wave propagation
based on an excitable reaction-diffusion equation, it has been found that the
normal velocity should be relevant to the local curvature \cite{KeS,Kur,TyK}.
For example, the so-called \textit{eikonal-curvature} relation or
\textit{mean-curvature flow} type equation (see (\ref{2.2}) below) has been
shown to be correct when linear approximation is considered \cite{KeS,Kur,TyK}.

We should note here that Wiener and Rosenblueth considered only pinned spirals
that rotate around an obstacle. However, subsequent numerical simulations and
experiments with various excitable media have revealed that the same media can
also support spiral waves free of any obstacles \cite{SkS,Win4}.

Later, Mikhailov \textit{et al.}\cite{MDZ} proposed the following standard
kinematic model:%

\begin{align}
\frac{\partial\kappa}{\partial t}+\frac{\partial\kappa}{\partial s}\left(
\int_{0}^{s}\kappa Vd\xi+G\right)   &  +\kappa^{2}V+\frac{\partial^{2}%
V}{\partial s^{2}}=0,\quad0\leq s<\infty,\quad0\leq t<\infty,\label{2.1}\\
V  &  =V_{0}-D\kappa,\quad0\leq s<\infty,\quad0\leq t<\infty,\label{2.2}\\
G  &  =G_{0}-r\kappa_{0},\text{ \ }0\leq t<\infty,\label{2.3}\\
\frac{d\kappa_{0}}{dt}  &  =\left.  -G\,\frac{\partial\kappa}{\partial
s}\right\vert _{s=0},\quad0<t<\infty,\quad\text{and}\label{2.4}\\
\lim_{s\rightarrow\infty}\kappa(s,t)  &  =0,\quad0<t<\infty, \label{2.5}%
\end{align}
where $\kappa=\kappa(s,t)$ is the curvature of a plane curve depending on time
$t$ and the arc length $s$ measured from the free tip, $V=V(s,t)$ is the
normal velocity, $G=G(t)$ is the tangential velocity of the tip, $\kappa
_{0}=\kappa_{0}(t):=\lim_{s\rightarrow0}\kappa(s,t),$ and $V_{0},$ $D,$
$G_{0},$ $r$ are parameters determined from the media.

Here, (\ref{2.1}) is a general equation satisfied by a plane curve of any
point at which the curve moves in the normal direction (see (4.41) in
\cite{Mer} or (2.7) in \cite{MDZ}); (\ref{2.2}) is the {eikonal-curvature}
relation, which relates the normal velocity of the curve to the local
curvature (see (4.10) in \cite{Mer} or (1.1) in \cite{MDZ}); (\ref{2.3}) is
the tangential velocity of the tip\footnote{(\ref{2.3}) is postulated in
\cite{BDM, BDZ} and derived by perturbation techniques in Appendix A of
\cite{MDZ}, when linear approximation is considered.} and (\ref{2.4}) and
(\ref{2.5}) are boundary conditions (see (2.10) in \cite{MDZ}).

We explain (\ref{2.3}) as follows. The curvature near the tip, hence
$\kappa_{0}$, affects the normal velocity near the tip by (\ref{2.2}). Then,
by (3.27) from \cite{Mer}, the normal velocity near the tip affects the width
of the excited region near the tip. Moreover, it has been observed that the
width of the excited region near the tip influences the tangential velocity of
the tip \cite{MGP,MPG}. Hence $\kappa_{0}$ is related to $G$.

Notice that (\ref{2.2}) and (\ref{2.3}) are good linear approximations of the
perturbation expansions only when the curvature is small. Hence ``$\kappa
\ll\upsilon_{0}$'' in (\ref{2.2}) (or $\upsilon\approx\upsilon_{0}$) is required.

\paragraph{Remarks}

Elkin \textit{et al.}\cite{EBH} used perturbations techniques and Fredholm
alternative theorem to obtain more general equations compared to (\ref{2.3})
and (\ref{2.4}) here (see (28), (32), (36), and (37) in \cite{EBH}). In fact,
it has been found that (\ref{2.3}) is only a special case of (32) in
\cite{EBH} (see the context below (37) in \cite{EBH}) and (\ref{2.4}) is not
correct. Later Elkin \textit{et al. }\cite{EBH2} solved the above new model
under some simplification, i.e., replacing $V$ with $V_{0}$ in the term
\textquotedblleft$\int_{0}^{s}\kappa Vd\xi$ \textquotedblright\ in (\ref{2.1})
(see (4) in \cite{EBH2}).

\bigskip

In this paper, we consider steadily rotating spiral waves (so $\frac
{\partial\kappa}{\partial t}=0$). Let $G$ be a constant and let $\omega$ be
the momentary rotational angular velocity of the tip. Then, using (\ref{2.2})
and integrating the both sides of (\ref{2.1}) give rise to:%

\begin{equation}
\kappa(s)\{\int_{0}^{s}k(\xi)[V_{0}-D\kappa(\xi)]d\xi+G\}-D\kappa^{\prime
}(s)=\omega. \label{2.6}%
\end{equation}
Thus it follows that
\begin{equation}
\kappa^{\prime}(0)=\frac{G}{D}\kappa(0)-\frac{\omega}{D}=\frac{G}{D}\kappa
_{0}-\frac{\omega}{D}, \label{2.7}%
\end{equation}
(see also (6) in \cite{EBH2}). After differentiating (\ref{2.6}) once and
replacing the integral term, we obtain:
\begin{equation}
-D\kappa^{\prime\prime}+\kappa^{2}(V_{0}-D\kappa)+\frac{\kappa^{\prime}%
}{\kappa}(\omega+D\kappa^{\prime})=0, \label{2.8}%
\end{equation}
provided $\kappa\neq0.$ Without loss of generality, we assume $\kappa>0$
hereafter, since if $\kappa<0,$ then (\ref{2.7}) and (\ref{2.8}) are invariant
by replacing $\kappa,$ $\omega,$ $V_{0}$ with $-\kappa,$ $-\omega,$ $-V_{0},$
respectively. Let $l(s):=\ln\kappa(s).$ Then (\ref{2.5}), (\ref{2.7}), and
(\ref{2.8}) give rise to:%
\begin{equation}
\frac{d^{2}l(s)}{ds^{2}}=-g(l(s))+\frac{\omega}{D}\frac{dl(s)}{ds}e^{-l(s)},
\label{2.9a}%
\end{equation}%
\begin{equation}
l(0)=\ln\kappa_{0}, \label{2.9b}%
\end{equation}%
\begin{equation}
l^{\prime}(0)=-\frac{\omega}{D\kappa_{0}}+\frac{G}{D}, \label{2.9c}%
\end{equation}%
\begin{equation}
\underset{s\rightarrow\infty}{\lim}l(s)=-\infty, \label{2.9d}%
\end{equation}
where $g(l):=e^{2l}-\frac{V_{0}}{D}e^{l}.$ The above equations (\ref{2.9a}),
(\ref{2.9b}), (\ref{2.9c}) are also equivalent to the following ones:%

\begin{equation}
l^{^{\prime}}(s)=v(s) \label{2.10a}%
\end{equation}%
\begin{equation}
v^{\prime}(s)=-g(l(s))+\frac{\omega}{D}e^{-l(s)}v(s) \label{2.10b}%
\end{equation}%
\begin{equation}
v(0)=-\frac{\omega}{D}e^{-l(0)}+\frac{G}{D}. \label{2.10c}%
\end{equation}
Furthermore, if $v(s)\neq0,$ then the solution $v(s)$ in (\ref{2.10a}%
),(\ref{2.10b}) can be viewed as a function of $l$ by the inverse function
theorem and satisfies the following equation:%

\begin{equation}
\frac{dv}{dl}=-\frac{g(l)}{v}+\frac{\omega}{D}e^{-l}. \label{2.11}%
\end{equation}

Note that in \cite{IIY}, Ishimura \textit{et al. }considered the above problem
(\ref{2.9a})-(\ref{2.9d}) with $G=0$ and $l^{\prime}(0)=0$ (see (13) in
\cite{IIY}). Therefore, the results in \cite{IIY} are valid only\ for the case
\textquotedblleft$\omega=0$\textquotedblright.

Now for any excitable media with given parameters $V_{0}>0,$ $D>0$ which are
fixed throughout this paper$,$ we want to find suitable $G,$ $\omega,$ and
$l(0)$ (i.e., $\kappa_{0}$) such that (\ref{2.9d}), (\ref{2.10a}),
(\ref{2.10b}), and (\ref{2.10c}) has global existence solutions$.$ Let
$E(l,v):=\frac{1}{2}v^{2}+\frac{1}{2}e^{2l}-\frac{V_{0}}{D}e^{l}.$

\section{Behavior of the Solutions}

\begin{lemma}
\label{Lemma0.1}For any $\omega<0,$ there are no global existence solutions of
(\ref{2.10a}),(\ref{2.10b}) satisfying $\underset{s\longrightarrow\infty}%
{\lim}l(s)=-\infty.$
\end{lemma}

\begin{proof}
Since $\omega<0,$ given any solution of (\ref{2.10a}),(\ref{2.10b}), say,
$(l(s;\omega),v(s;\omega))$ (dropping out the dependence on initial data here
for simplicity)$,$ $E(l(\cdot;\omega),v(\cdot;\omega))$ is a decreasing
function by (11) of [1]. Note also that $E(l,v)\geq E(\ln\frac{V_{0}}%
{D},0)=\frac{-1}{2}(\frac{V_{0}}{D})^{2}.$ Thus $E(l(\cdot;\omega
),v(\cdot;\omega))$ is bounded. If $l(\cdot;\omega)$ exists globally with
$\underset{s\longrightarrow\infty}{\lim}l(s;\omega)=-\infty$, then we have
$\underset{s\rightarrow\infty}{\lim}se^{2l(s)}=\frac{\omega}{2V_{0}}$ by LEMMA
3 of \cite{IIY}. This is impossible, since $\underset{s\rightarrow\infty}%
{\lim}se^{2l(s)}\geq0>\frac{\omega}{2V_{0}}.$ This completes the proof.
\end{proof}

\bigskip

By above lemma and THEOREM 2 of \cite{IIY}, we only need to consider
$\omega>0$ now.

\begin{lemma}
\label{Lemma0.3}Given any global existence solution of (\ref{2.10a}%
),(\ref{2.10b}), $(l(s;\omega),v(s;\omega)),$ if $v(0;\omega)>0,$ then there
is some $\sigma>0$ such that $v(\sigma;\omega)=0$ (no matter the sign of
$\omega$)$.$
\end{lemma}

\begin{proof}
Suppose not, i.e., $(l(s;\omega),v(s;\omega))$ always stays in the upper half
plane: $v>0.$ Note that the upper half plane is splitted into three parts:
$R_{1}:=\{(l,v):v>0$ and $-g(l)+\frac{\omega}{D}e^{-l}v>0\},$ $R_{2}%
:=\{(l,v):v>0$ and $-g(l)+\frac{\omega}{D}e^{-l}v=0\},$ and $R_{3}%
:=\{(l,v):v>0$ and $-g(l)+\frac{\omega}{D}e^{-l}v<0\}.$ If $(l(s;\omega
),v(s;\omega))$ always stays in $R_{1},$ then both $l(\cdot;\omega
),v(\cdot;\omega)$ increase as $s\rightarrow\infty$ by (\ref{2.10a}%
),(\ref{2.10b}) and so both $\underset{s\rightarrow\infty}{\lim}$
$l(s;\omega)$ and\textbf{ }$\mu:=\underset{s\rightarrow\infty}{\lim}$
$v(s;\omega)>0$ exist. Thus clearly, we have $\underset{s\rightarrow\infty
}{\lim}$ $l(s;\omega)=\infty.$ Since $v^{^{\prime}}(s;\omega)>0$ and
$\underset{s\rightarrow\infty}{\lim}$ $-g(l(s;\omega))=-\infty$ now, we have
$\mu=\infty$ by (\ref{2.10b}). On the other hand, it follows from
(\ref{2.10b}) that:%
\begin{equation}
v(s;\omega)+\frac{\omega}{D}e^{-l(s;\omega)}-(v(0;\omega)+\frac{\omega}%
{D}e^{-l(0;\omega)})=-\int_{0}^{s}g(l(t;\omega))dt. \label{3.1}%
\end{equation}
Taking limit as $s\rightarrow\infty$ on the both side of (\ref{3.1}), we have%
\begin{equation}
\underset{s\rightarrow\infty}{\lim}[v(s;\omega)+\frac{\omega}{D}%
e^{-l(s;\omega)}-(v(0;\omega)+\frac{\omega}{D}e^{-l(0;\omega)})]=\infty
=-\int_{0}^{\infty}g(l(t;\omega))dt.\text{ } \label{3.2}%
\end{equation}
However, $-\int_{0}^{\infty}g(l(t;\omega))dt=-\int_{0}^{\delta}g(l(t;\omega
))dt-\int_{\delta}^{\infty}g(l(t;\omega))dt,$ where $l(\delta;\omega)=\ln
\frac{V_{0}}{D}.$ Since $g(l)\geq0$ as $l\geq\ln\frac{V_{0}}{D}$ and
$l(t;\omega)$ $\nearrow\infty$ as $t\nearrow\infty$, we have $-\int_{\delta
}^{\infty}g(l(t;\omega))dt<0$ and so
\begin{equation}
-\int_{0}^{\infty}g(l(t;\omega))dt<-\int_{0}^{\delta}g(l(t;\omega))dt<\infty.
\label{3.3}%
\end{equation}
Thus (\ref{3.3}) contradicts (\ref{3.2}). Now the only possible case is that
$(l(s;\omega),v(s;\omega))$ goes through the curve $R_{2}$, enters the region
$R_{3},$ and then stays there for ever. Hence $l(\cdot;\omega)$ and
$v(\cdot;\omega)$ are increasing and decreasing functions respectively now by
(\ref{2.10a}),(\ref{2.10b}). Clearly $\mu:=\underset{s\rightarrow\infty}{\lim
}$ $v(s;\omega)\in\lbrack0,\infty),$ since $v(s;\omega)>0.$ Note that the
\textquotedblleft Poincare-Bendixson theorem\textquotedblright\ implies
$\underset{s\rightarrow\infty}{\lim}$ $l(s;\omega)=\infty;$ otherwise,
$(l(s;\omega),v(s;\omega))$ will stay in a bounded region and so
$\underset{s\rightarrow\infty}{\lim}$ $l(s;\omega)=$ $\ln\frac{V_{0}}{D}$ and
$\underset{s\rightarrow\infty}{\lim}$ $v(s;\omega)=0,$ since $(\ln\frac{V_{0}%
}{D},0)$ is the only possible $\omega-$limit point of $(l(s;\omega
),v(s;\omega)).$ A contradiction then occurs, since $\underset{s\rightarrow
\infty}{\lim}$ $l(s;\omega)\neq\ln\frac{V_{0}}{D}$ by (\ref{2.10a}). Now it
follows from \textquotedblleft$\underset{s\rightarrow\infty}{\lim}$
$l(s;\omega)=\infty$\textquotedblright\ and (\ref{2.10b}) that $\underset
{s\rightarrow\infty}{\lim}v^{^{\prime}}(s;\omega)=-\infty$ which contradicts
\textquotedblleft$v(s;\omega)>0,\forall s\geq0$\textquotedblright.
\end{proof}

\bigskip

\begin{lemma}
\label{Lemma0.2}For $\omega>0,$ given any solution of (\ref{2.10a}%
),(\ref{2.10b}), $(l(s;\omega),v(s;\omega)),$ with maximal existence interval
$[0,\Gamma),$ if $E(l(s^{\#};\omega),v(s^{\#};\omega))\geq0$ for some
$s^{\#}\in\lbrack0,\Gamma),$ then $\Gamma<\infty$ and $\underset
{s\rightarrow\Gamma^{-}}{\lim}$ $l(s;\omega)=\underset{s\rightarrow\Gamma^{-}%
}{\lim}v(s;\omega)=-\infty,$ i.e., blows up in finite $s.$
\end{lemma}

\begin{proof}
Assume the contrary that $\Gamma=\infty.$ Then by Lemma \ref{Lemma0.3},
$(l(s;\omega),v(s;\omega))$ will enter the lower half plane after a short
time, say, $\sigma$. Since $E(l(\cdot;\omega),v(\cdot;\omega))$ is increasing
and $E(l(s^{\#};\omega),v(s^{\#};\omega))\geq0$ for some $s^{\#}\in
\lbrack0,\Gamma),$ we have $E(l(s;\omega),v(s;\omega))>0$ $\forall s\in
(s^{\#},\infty).$ It then follows that $v(s;\omega)<0,$ $\forall s\in
(\max\{\sigma,s^{\#}\},\infty),$ since the level curve $E(l,v)=0$ is a barrier
preventing $(l(s;\omega),v(s;\omega))$ from touching the $l-$axis. Clearly,
both $\mu_{1}:=\underset{s\rightarrow\infty}{\lim}l(s;\omega)$ $<\infty$ and
$\mu_{2}:=\underset{s\rightarrow\infty}{\lim}v(s;\omega)\leq0$ exists. Note
that $\mu_{2}\neq0;$ otherwise, again by the \textquotedblleft
Poincare-Bendixson theorem\textquotedblright\ we have $\mu_{1}=-\infty$ and so
$\underset{s\rightarrow\infty}{\lim}$ $E(l(s;\omega),v(s;\omega))=0$ which is
impossible. Now it follows from \textquotedblleft$\mu_{2}\in\lbrack
-\infty,0)"$ that $\mu_{1}=-\infty.$ By (\ref{2.10b}), we then have
$\underset{s\rightarrow\infty}{\lim}v^{^{\prime}}(s;\omega)=-\infty$ and so
$\mu_{2}=-\infty.$ Moreover, we also have $\underset{s\rightarrow\infty}{\lim
}v^{^{^{\prime\prime}}}(s;\omega)=-\infty$ by differentiating both sides of
(\ref{2.10b}) and taking limit.

Multiplying both sides of (\ref{3.1}) by $e^{l(s;\omega)}$ and integrating
them$,$ we then have%
\begin{equation}
\frac{\omega s}{D}+\int_{l(0;\omega)}^{l(s;\omega)}e^{l}dl-(v(0;\omega
)+\frac{\omega}{D}e^{-l(0;\omega)})\int_{0}^{s}e^{l(p;\omega)}dp=-\int_{0}%
^{s}e^{l(p;\omega)}\int_{0}^{p}g(l(t;\omega))dtdp. \label{3.4}%
\end{equation}
Now we want to show that $e^{l(p;\omega)}=o(\frac{1}{p^{3}})$ as
$p\rightarrow\infty.$ Utilizing l'H\^{o}pital's rule, we have%
\begin{align}
\underset{p\rightarrow\infty}{\lim}\frac{e^{l(p;\omega)}}{p^{-3}}  &
=\underset{p\rightarrow\infty}{\lim}\frac{p^{3}}{e^{-l(p;\omega)}}%
=\underset{p\rightarrow\infty}{\lim}\frac{3p^{2}}{-v(p;\omega)e^{-l(p;\omega
)}}=\underset{p\rightarrow\infty}{\lim}\frac{6p}{-v^{^{\prime}}(p;\omega
)e^{-l(p;\omega)}+v^{2}(p;\omega)e^{-l(p;\omega)}}\label{3.5}\\
&  =\underset{p\rightarrow\infty}{\lim}\frac{6}{[-v^{^{\prime\prime}}%
(p;\omega)+3v^{^{\prime}}(p;\omega)v(p;\omega)-v^{3}(p;\omega)]e^{-l(p;\omega
)}}=0.\nonumber
\end{align}
Since $g(l)\geq g(\ln\frac{V_{0}}{2D}),$ we have $-\int_{0}^{\infty
}e^{l(p;\omega)}\int_{0}^{p}g(l(t;\omega))dtdp\leq-g(\ln\frac{V_{0}}{2D}%
)\int_{0}^{\infty}e^{l(p;\omega)}pdp<\infty$ by (\ref{3.5}). Also,
$\underset{s\rightarrow\infty}{\lim}\int_{l(0;\omega)}^{l(s;\omega)}%
e^{l}dl=\int_{l(0;\omega)}^{-\infty}e^{l}dl\in(-\infty,\infty)$ and $\int
_{0}^{\infty}e^{l(p;\omega)}dp\in(0,\infty)$ by (\ref{3.5}) again. By taking
limit as $s\rightarrow\infty$ on the both sides of (\ref{3.4})$,$ a
contradictions occurs. Therefore $\Gamma<\infty$ and so the proof is finished.
\end{proof}

\bigskip

\begin{lemma}
\label{Lemma1}Given any $\omega\in(0,\frac{2V_{0}^{2}}{D}),$ there exists a
unique $l^{\ast}=$ $l^{\ast}(\omega)\in(\ln\frac{V_{0}}{D},\ln\frac{2V_{0}}%
{D})$ such that (\ref{2.10a}),(\ref{2.10b}) has a global existence solution
$(\widetilde{l}(s;\omega),\widetilde{v}(s;\omega))$ defined on $(-\infty
,\infty)$ with $(\widetilde{l}(0;\omega),\widetilde{v}(0;\omega))=(l^{\ast
},0),$ $\widetilde{l}(s;\omega)\searrow-\infty,$ $\widetilde{v}(s;\omega
)\rightarrow0$ as $s\longrightarrow\infty,$ and $\widetilde{v}(s;\omega)<0$
for all $s>0$. Moreover, $(\widetilde{l}(s;\omega),\widetilde{v}(s;\omega))$
rotates counterclockwise around $(\ln\frac{V_{0}}{D},0)$ and tends to it as
$s$ decreases to $-\infty.$
\end{lemma}

\begin{proof}
By THEOREM 3, LEMMA 4, and LEMMA 5 in [1], we just know there exists an
interval $[\sup A,\inf B]\subset(\ln\frac{V_{0}}{D},\ln\frac{2V_{0}}{D}]$ such
that the solution of (\ref{2.10a}),(\ref{2.10b}) with initial data
$(l_{0},0),$ where $l_{0}\in$ $[\sup A,\inf B],$ has the properties mentioned
in this lemma, denoted by $(\widetilde{l}(s;\omega,l_{0}),\widetilde
{v}(s;\omega,l_{0}))$. Note that by Lemma \ref{Lemma0.2}, we have $[\sup
A,\inf B]\subset(\ln\frac{V_{0}}{D},\ln\frac{2V_{0}}{D}).$ Now we want to
prove $\sup A=\inf B.$

Since the solution $\widetilde{v}(s;\omega,l_{0})\neq0$, for all
$s\in(0,\infty),$ $\widetilde{v}(s;\omega,l_{0})$ can be viewed as a function
of $l$ by the inverse function theorem, i.e., $\widetilde{v}(s;\omega
,l_{0})=\widetilde{v}(l;\omega,l_{0}),$ for $l\in\lbrack-\infty,l_{0}]$ and so
satisfies (\ref{2.11}). For any two solutions $\widetilde{v}(l;\omega
,l_{01}),$ $\widetilde{v}(l;\omega,l_{02}),$ denoted by $\widetilde{v}_{1}(l)$
and $\widetilde{v}_{2}(l)$ respectively for simplicity, by (\ref{2.11}) we
have
\begin{equation}
\frac{d(\widetilde{v}_{1}-\widetilde{v}_{2})}{dl}=\frac{g(l)}{\widetilde
{v}_{1}\widetilde{v}_{2}}(\widetilde{v}_{1}-\widetilde{v}_{2}). \label{3.6}%
\end{equation}
Thus $(\widetilde{v}_{1}-\widetilde{v}_{2})(l)=$ $(\widetilde{v}%
_{1}-\widetilde{v}_{2})(\ln\frac{V_{0}}{D})\exp(-\int_{l}^{\ln\frac{V_{0}}{D}%
}\frac{g(l)}{\widetilde{v}_{1}(l)\widetilde{v}_{2}(l)}dl)$ for $l<\ln
\frac{V_{0}}{D}.$ Suppose ($\widetilde{v}_{1}-\widetilde{v}_{2})(\ln
\frac{V_{0}}{D})\neq0.$ Then since $\widetilde{v}_{1},\widetilde{v}_{2}<0$ and
$g(l)<0$ for $l<\ln\frac{V_{0}}{D},$ we have$\underset{l\longrightarrow
-\infty}{\lim}(\widetilde{v}_{1}-\widetilde{v}_{2})(l)\neq0.$This is
impossible, since $\underset{s\longrightarrow\infty}{\lim}\widetilde
{v}(s;\omega,l_{01})=\underset{s\longrightarrow\infty}{\lim}\widetilde
{v}(s;\omega,l_{02})=0.$ By uniqueness, we have $l_{01}=l_{02}$ and so $\sup
A=\inf B.$ Let $l^{\ast}(\omega):=\sup A=\inf B$ and $(\widetilde{l}%
(s;\omega),\widetilde{v}(s;\omega)):=(\widetilde{l}(s;\omega,l^{\ast}%
(\omega)),\widetilde{v}(s;\omega,l^{\ast}(\omega))).$

Finally, since the $\alpha-$limit set of $(\widetilde{l}(s;\omega
),\widetilde{v}(s;\omega))$ consists of only one point $(\ln\frac{V_{0}}%
{D},0)$ which is a unstable focus, $(\widetilde{l}(s;\omega),\widetilde
{v}(s;\omega))$ is also defined on $(-\infty,0]$ and $\underset{s\rightarrow
-\infty}{\lim}(\widetilde{l}(s;\omega),\widetilde{v}(s;\omega))=(\ln
\frac{V_{0}}{D},0).$
\end{proof}

\bigskip

Given any solution of (\ref{2.10a}),(\ref{2.10b}), say, $(l(s;\omega
),v(s;\omega)),$ $E(l(\cdot;\omega),v(\cdot;\omega))$ is an increasing
function, since $\omega>0$. Also, by LEMMA 5 in \cite{IIY}, $(\widetilde
{l}(s;\omega),\widetilde{v}(s;\omega))$ always stays inside the curve:
$E(l,v)=0$ in the $l-v$ phase plane$.$

\begin{proposition}
\label{Proposition2}Suppose $\omega\in(0,\frac{2V_{0}^{2}}{D}).$ For any
global existence solution of (\ref{2.10a}),(\ref{2.10b}), its orbit curve in
the $l-v$ plane coincides with that of $(\widetilde{l}(s;\omega),\widetilde
{v}(s;\omega))$, i.e., they both lie on the same solution curve in the $l-v$ plane.
\end{proposition}

\begin{proof}
Let $(l(s;\omega),v(s;\omega))$ be a global existence solution of
(\ref{2.10a}),(\ref{2.10b}), i.e.,existence interval contains $[0,\infty),$
(also dropping out the dependence on initial data here for simplicity).
Suppose $(l(s;\omega),v(s;\omega))$ stays in a bounded region in the $l-v$
plane. Since the equilibrium $(\ln a,0)$ is a unstable focus, the only
possible $\omega-$limit set of $(l(s;\omega),v(s;\omega))$ consists of
periodic orbits by the \textquotedblleft Poincare-Bendixson
theorem\textquotedblright. However, by the \textquotedblleft Bendixson
criterion\textquotedblright\ (see [ ]), there is no any periodic solution of
(\ref{2.10a}),(\ref{2.10b}). Thus our assumption is wrong. Therefore, any
global existence solution is unbounded.

Note that by Lemma \ref{Lemma0.2}, $(l(s;\omega),v(s;\omega))$ always stays
inside the curve:$E(l,v)=0$ in the $l-v$ phase plane$.$ Since $(l(s;\omega
),v(s;\omega))$ is unbounded, we can choose some $s_{1}$ large enough such
that $(l(s_{1};\omega),v(s_{1};\omega))$ is located in the lower
\textquotedblleft left tail\textquotedblright\ region : $\{(l,v)\in
R^{2}\left\vert E(l,v)<0,\text{ }v<0,\text{ and }l<l_{L},\text{ for some small
enough }l_{L}\right.  \}$ with $E(l(s_{1};\omega),v(s_{1};\omega))>E(l^{\ast
},0)$ by the structure of level curve of energy $E(l,v)$ with negative value,
\textit{i.e., }$E(l,v)=c<0.$ Suppose $v(s_{2};\omega)=0$ for some $s_{2}%
>s_{1}.$ Then by the unboundedness of $(l(s;\omega),v(s;\omega)),$ the
direction of flow on the upper half plane, and the increasingness of
$E(l(\cdot;\omega),v(\cdot;\omega)),$ there exists some $s_{3}>s_{2}$ such
that $v(s_{3};\omega)=0$ and $l(s_{3};\omega)$ $\in(l^{\ast},\ln\frac{2V_{0}%
}{D}).$ Note that any solution curve which intersects the $l-$axis at
$(l_{R},0)$ with $l_{R}\in(l^{\ast},\ln\frac{2V_{0}}{D})$ must touch the level
curve $E(l,v)=0$ and so blows up in finite time by Lemma \ref{Lemma0.2}, since
$l^{\ast}=\inf B$ (see \cite{IIY} for the definition of set $B$) and the
uniqueness theorem of ordinary differential equation. Thus $(l(s;\omega
),v(s;\omega))$ blows up in finite time which leads to a contradiction.
Therefore our assumption is wrong. We conclude that $v(s;\omega)<0$ for all
$s\geq s_{1}$ which leads to the following facts: $\lim_{s\longrightarrow
\infty}l(s;\omega)=-\infty$ and $\lim_{s\longrightarrow\infty}v(s;\omega)=0.$

Now using similar argument as in Lemma \ref{Lemma1}, the solution curve of
$(l(s;\omega),v(s;\omega))$ must coincide with that of $(\widetilde
{l}(s;\omega),\widetilde{v}(s;\omega))$ in the $l-v$ plane.
\end{proof}

\bigskip

Let $l_{1R}^{\ast}(\omega):=l^{\ast}(\omega),$ $l_{2R}^{\ast}(\omega
),...,l_{iR}^{\ast}(\omega),...$ $>\ln\frac{V_{0}}{D}$ and $l_{1L}^{\ast
}(\omega),$ $l_{2L}^{\ast}(\omega),...,l_{iL}^{\ast}(\omega),...$ $<\ln
\frac{V_{0}}{D}$ be the sequence of intersection points of $(\widetilde
{l}(s;\omega),\widetilde{v}(s;\omega))$ with the $l-$axis greater than and
less than $\ln\frac{V_{0}}{D}$ respectively as $s\rightarrow-\infty.$ Note
that $l_{iR}^{\ast}(\omega)>l_{jR}^{\ast}(\omega)$ and $l_{iL}^{\ast}%
(\omega)<l_{jL}^{\ast}(\omega)$ as $i<j.$

\begin{lemma}
\begin{enumerate}
\item[(i).] \label{Lemma4}$l_{iR}^{\ast}(\cdot)$ and $l_{iL}^{\ast}(\cdot)$
are decreasing and increasing functions respectively on $(0,\frac{2V_{0}^{2}%
}{D}).$

\item[(ii).] $l_{iR}^{\ast}(\cdot)$ and $l_{iL}^{\ast}(\cdot)$ are continuous
functions on $(0,\frac{2V_{0}^{2}}{D})$.

\item[(iii).] $\underset{\omega\rightarrow0^{+}}{\lim}l_{iR}^{\ast}%
(\omega)=\ln\frac{2V_{0}}{D}$ and $\underset{\omega\rightarrow0^{+}}{\lim
}l_{iL}^{\ast}(\omega)=-\infty,$ $i=1,2,3,...$
\end{enumerate}
\end{lemma}

\begin{proof}
Since $\widetilde{v}(s;\omega_{1})\neq0,\widetilde{v}(s;\omega_{2})\neq0$, for
all $s\in(0,\infty),$ $\widetilde{v}(s;\omega_{1}),\widetilde{v}(s;\omega
_{2})$ can be viewed as a function of $l$ by the inverse function theorem,
i.e., $\widetilde{v}(s;\omega_{1})=\widetilde{v}(l;\omega_{1}),\widetilde
{v}(s;\omega_{2})=\widetilde{v}(l;\omega_{2})$ for $l\in\lbrack-\infty
,l_{1R}^{\ast}(\omega_{1})),l\in\lbrack-\infty,l_{1R}^{\ast}(\omega_{2}))$
respectively. Note that for $\omega_{1},\omega_{2},$ the corresponding
solutions $v_{1}(l):=v(l;\omega_{1}),$ $v_{2}(l):=v(l;\omega_{2})$ of
(\ref{2.11}) have the following relation:%
\[
(v_{1}-v_{2})(l)=
\]%
\begin{equation}
(v_{1}-v_{2})(\tau)\exp(\int_{\tau}^{l}\frac{g(t)}{v_{1}(t)v_{2}(t)}%
dt)+(\frac{\omega_{1}-\omega_{2}}{D})\exp(\int_{\tau}^{l}\frac{g(t)}%
{v_{1}(t)v_{2}(t)}dt)\int_{\tau}^{l}\exp(-t-\int_{\tau}^{t}\frac{g(\rho
)}{v_{1}(\rho)v_{2}(\rho)}d\rho)dt, \label{3.7}%
\end{equation}
whenever both $v_{1},v_{2}$ exist on $[\tau,l]$ or$\ [l,\tau].$ Now apply
(\ref{3.7}) to $\widetilde{v}_{1}(l):=\widetilde{v}(l;\omega_{1}),$
$\widetilde{v}_{2}(l):=\widetilde{v}(l;\omega_{2}),$ for $\omega_{1}%
,\omega_{2}\in$ $(0,\frac{2V_{0}^{2}}{D})$ with $\omega_{1}<\omega_{2}.$

First we consider $l_{1R}^{\ast}(\cdot).$

\textbf{Case 1. }$l_{1R}^{\ast}(\omega_{1})<l_{1R}^{\ast}(\omega_{2}).$

Then there is a small $\varepsilon>0$ such that $\widetilde{v}_{1}%
(l),\widetilde{v}_{2}(l)$ are defined on $(-\infty,l_{1R}^{\ast}(\omega
_{1})-\varepsilon]$ with $\widetilde{v}_{1}(l_{1R}^{\ast}(\omega
_{1})-\varepsilon)>\widetilde{v}_{2}(l_{1R}^{\ast}(\omega_{1})-\varepsilon).$
Taking $\tau=l_{1R}^{\ast}(\omega_{1})-\varepsilon$ and $l\leq\tau,$ by
(\ref{3.7}) we obtain:
\[
(\widetilde{v}_{1}-\widetilde{v}_{2})(l)\geq(\widetilde{v}_{1}-\widetilde
{v}_{2})(\tau)\exp(\int_{\tau}^{l}\frac{g(t)}{\widetilde{v}_{1}(t)\widetilde
{v}_{2}(t)}dt).
\]
Since $\underset{l\rightarrow-\infty}{\lim}(\widetilde{v}_{1}-\widetilde
{v}_{2})(\tau)\exp(\int_{\tau}^{l}\frac{g(t)}{\widetilde{v}_{1}(t)\widetilde
{v}_{2}(t)}dt)>0,$ we have $\underset{l\rightarrow-\infty}{\lim}(\widetilde
{v}_{1}-\widetilde{v}_{2})(l)>0.$ A contradiction occurs.

\textbf{Case 2. }$l_{1R}^{\ast}(\omega_{1})=l_{1R}^{\ast}(\omega_{2})$ and
there is $\epsilon>0$ such that $\widetilde{v}_{1}(l_{1R}^{\ast}(\omega
_{1})-\epsilon)\geq\widetilde{v}_{2}(l_{1R}^{\ast}(\omega_{1})-\epsilon).$

Then taking $\tau=l_{1R}^{\ast}(\omega_{1})-\epsilon$ and $l\leq\tau,$ by
(\ref{3.7}) we obtain:%
\[
(\widetilde{v}_{1}-\widetilde{v}_{2})(l)\geq(\frac{\omega_{1}-\omega_{2}}%
{D})\exp(\int_{\tau}^{l}\frac{g(t)}{\widetilde{v}_{1}(t)\widetilde{v}_{2}%
(t)}dt)\int_{\tau}^{l}\exp(-t-\int_{\tau}^{t}\frac{g(\rho)}{\widetilde{v}%
_{1}(\rho)\widetilde{v}_{2}(\rho)}d\rho)dt.
\]
Since $\underset{l\rightarrow-\infty}{\lim}(\frac{\omega_{1}-\omega_{2}}%
{D})\exp(\int_{\tau}^{l}\frac{g(t)}{\widetilde{v}_{1}(t)\widetilde{v}_{2}%
(t)}dt)\int_{\tau}^{l}\exp(-t-\int_{\tau}^{t}\frac{g(\rho)}{\widetilde{v}%
_{1}(\rho)\widetilde{v}_{2}(\rho)}d\rho)dt>0,$ we also have a contradiction.
Hence the only possible case is $l_{1R}^{\ast}(\omega_{1})=l_{1R}^{\ast
}(\omega_{2})$ and $\widetilde{v}_{1}(l)<\widetilde{v}_{2}(l)$ for
$l\in(-\infty,l_{1R}^{\ast}(\omega_{1}))$ other than $l_{1R}^{\ast}(\omega
_{1})>l_{1R}^{\ast}(\omega_{2}).$ Taking any $\tau\in\lbrack\ln\frac{V_{0}}%
{D},l_{1R}^{\ast}(\omega_{1}))$ and $l\in\lbrack\tau,l_{1R}^{\ast}(\omega
_{1})),$ by (\ref{3.7}) we have
\begin{equation}
(\widetilde{v}_{1}-\widetilde{v}_{2})(l)<(\frac{\omega_{1}-\omega_{2}}{D}%
)\exp(\int_{\tau}^{l}\frac{g(t)}{\widetilde{v}_{1}(t)\widetilde{v}_{2}%
(t)}dt)\int_{\tau}^{l}\exp(-t-\int_{\tau}^{t}\frac{g(\rho)}{\widetilde{v}%
_{1}(\rho)\widetilde{v}_{2}(\rho)}d\rho)dt<0.\text{ } \label{3.8}%
\end{equation}
Clearly, the right hand side of (\ref{3.8}) will not tend to $0,$ as
$l\rightarrow l_{1R}^{\ast}(\omega_{1})^{-}.$ This contradicts to
$\underset{l\rightarrow l_{1R}^{\ast}(\omega_{1})^{-}}{\lim}(\widetilde{v}%
_{1}-\widetilde{v}_{2})(l)=0,$ since $l_{1R}^{\ast}(\omega_{1})=l_{1R}^{\ast
}(\omega_{2}).$

By above discussions, we conclude that $l_{1R}^{\ast}(\omega_{1})>l_{1R}%
^{\ast}(\omega_{2}),$ if \ $\omega_{1}<\omega_{2}.$

Now we consider $l_{1L}^{\ast}(\cdot).$ Given any $s_{L_{1}}<0$ such that the
solution curves of $(\widetilde{l}(s;\omega_{1}),\widetilde{v}(s;\omega_{1}))$
and $(\widetilde{l}(s;\omega_{2}),\widetilde{v}(s;\omega_{2}))$ stay in the
upper half plane when $s\in\lbrack s_{L_{1}},0),$ both $\widetilde{v}%
(s;\omega_{1})$ and $\widetilde{v}(s;\omega_{2})$ then can also be viewed as
functions of $l,$ $\widetilde{v}_{1}(\cdot)$, $\widetilde{v}_{2}(\cdot)$ as
before and so (\ref{3.7}) is also valid. Due to $l_{1R}^{\ast}(\omega
_{1})>l_{1R}^{\ast}(\omega_{2}),$ we can take $\tau=l_{1R}^{\ast}(\omega
_{2})-\delta$ for some small $\delta>0$ such that both $\widetilde{v}%
_{1}(\cdot)$ and $\widetilde{v}_{2}(\cdot)$ are defined on $(\max
\{l_{1L}^{\ast}(\omega_{1}),l_{1L}^{\ast}(\omega_{2})\},\tau]$ with
$\widetilde{v}_{1}(\tau)>$ $\widetilde{v}_{2}(\tau).$ By (\ref{3.7}) we have
\begin{equation}
(\widetilde{v}_{1}-\widetilde{v}_{2})(l)>(\frac{\omega_{1}-\omega_{2}}{D}%
)\exp(\int_{\tau}^{l}\frac{g(t)}{\widetilde{v}_{1}(t)\widetilde{v}_{2}%
(t)}dt)\int_{\tau}^{l}\exp(-t-\int_{\tau}^{t}\frac{g(\rho)}{\widetilde{v}%
_{1}(\rho)\widetilde{v}_{2}(\rho)}d\rho)dt>0, \label{3.9}%
\end{equation}
for $l\in$ $(\max\{l_{1L}^{\ast}(\omega_{1}),l_{1L}^{\ast}(\omega_{2}%
)\},\tau).$ Since the right hand side of (\ref{3.9}) will not tend to $0$ as
$l\rightarrow\max\{l_{1L}^{\ast}(\omega_{1}),l_{1L}^{\ast}(\omega_{2})\}^{+},$
we conclude that $l_{1L}^{\ast}(\omega_{1})<l_{1L}^{\ast}(\omega_{2}),$ if
\ $\omega_{1}<\omega_{2}.$

Arguing in the same way as above, we can derive $l_{iR}^{\ast}(\omega
_{1})>l_{iR}^{\ast}(\omega_{2})$ from ``$l_{(i-1)L}^{\ast}(\omega
_{1})<l_{(i-1)L}^{\ast}(\omega_{2})$''\ and then $l_{iL}^{\ast}(\omega
_{1})<l_{iL}^{\ast}(\omega_{2})$ from ``$l_{iR}^{\ast}(\omega_{1}%
)>l_{iR}^{\ast}(\omega_{2})$''\ for $i=2,3,4,...$ This completes (i).

Now we want to prove (ii). First consider $l_{1R}^{\ast}(\cdot).$ Note that by
(i) and Lemma \ref{Lemma1} $l_{0}:=\underset{\omega\rightarrow\omega_{0}}%
{\lim}$ $l_{1R}^{\ast}(\omega)$ exists for any given $\omega_{0}\in
(0,\frac{2V_{0}^{2}}{D})$ and $l_{0}\in$ $(\ln\frac{V_{0}}{D},\ln\frac{2V_{0}%
}{D}]$. Suppose $l_{0}=\ln\frac{2V_{0}}{D},$ i.e., the initial data and
parameter pair $((l_{1R}^{\ast}(\omega),0),\omega)$ tends to $((\ln
\frac{2V_{0}}{D},0),\omega_{0}).$ Then by the \textquotedblleft continuous
dependence on the initial data and parameter\textquotedblright\ theorem of
ordinary differential equation and Lemma \ref{Lemma0.2}, a contradiction
occurs, since $E(\ln\frac{2V_{0}}{D},0)=0$ and $E(l(\cdot;\omega
),v(\cdot;\omega))$ is increasing. Therefore we have $l_{0}\in(\ln\frac{V_{0}%
}{D},\ln\frac{2V_{0}}{D}).$ Also by the \textquotedblleft continuous
dependence\textquotedblright\ property, the solution of (\ref{2.10a}%
),(\ref{2.10b}) with initial data and parameter pair $((l_{0},0),\omega_{0})$
must stay in the region $\left\{  (l,v)\in R^{2}\left\vert E(l,v)<0\text{ and
}v<0\right.  \right\}  .$ Then by LEMMA 1 in [1] this solution is global
existence one. Thus by Lemma \ref{Lemma1} and Proposition \ref{Proposition2},
we have $l_{0}=l_{1R}^{\ast}(\omega_{0})$ and so $\underset{\omega
\rightarrow\omega_{0}}{\lim}$ $l_{1R}^{\ast}(\omega)=l_{1R}^{\ast}(\omega
_{0}).$ Hence $l^{\ast}(\cdot)=$ $l_{1R}^{\ast}(\cdot)$ is a continuous
function on $(0,\frac{2V_{0}^{2}}{D}).$

For $l_{iL}^{\ast}(\cdot),i=1,2,3,...$and $l_{iR}^{\ast}(\cdot),i=2,3,...$, by
the continuity of $l_{1R}^{\ast}(\cdot)$ and the \textquotedblleft continuous
dependence\textquotedblright\ property of o.d.e. when choosing long enough
existence intervals $[s_{L_{2i}},0],[s_{L_{2i-1}},0]$ such that the orbit
curve of $(\widetilde{l}(s;\omega_{0}),\widetilde{v}(s;\omega_{0}))$
intersects the $l-$axis $2i$ and $2i-1$ times respectively, both $l_{iL}%
^{\ast}(\omega)$ and $l_{iR}^{\ast}(\omega)$ will stay close to $l_{iL}^{\ast
}(\omega_{0})$ and $l_{iR}^{\ast}(\omega_{0})$ respectively as $\omega$ is
close to $\omega_{0}$. Therefore $l_{iL}^{\ast}(\cdot)$ and $l_{iR}^{\ast
}(\cdot)$ are continuous functions.

Finally we prove (iii). Samely, the limit $\zeta:=\underset{\omega
\rightarrow0^{+}}{\lim}l_{iR}^{\ast}(\omega)$ exists with $\zeta\in(\ln
\frac{V_{0}}{D},\ln\frac{2V_{0}}{D}].$ Suppose $\zeta\in(\ln\frac{V_{0}}%
{D},\ln\frac{2V_{0}}{D}).$ Consider the solution of (\ref{2.10a}%
),(\ref{2.10b}) with initial value-parameter pair $((\zeta,0),0)$; it is
periodic with some period, say, $T\in(0,\infty).$ Applying the
\textquotedblleft continuous dependence\textquotedblright\ \ property of
o.d.e. when choosing existence interval $[0,iT]$, solution $(\widetilde
{l}^{\#}(s;\omega),\widetilde{v}^{\#}(s;\omega)):=$ $(\widetilde
{l}(s+s_{i(\omega)};\omega),\widetilde{v}(s+s_{i(\omega)};\omega))$ over the
existence interval $[0,iT]$ where $(\widetilde{l}(s_{i(\omega)};\omega
),\widetilde{v}(s_{i(\omega)};\omega))=(l_{iR}^{\ast}(\omega),0),$ must
intersect the \textquotedblleft negative\textquotedblright\ $l-$axis (less
than $\ln\frac{V_{0}}{D}$) at least $i$ times as $\omega\rightarrow0^{+}$.
This is a contradiction. Thus $\underset{\omega\rightarrow0^{+}}{\lim}%
l_{iR}^{\ast}(\omega)=\ln\frac{2V_{0}}{D}$.

By (i) again, we have $\varsigma:=\underset{\omega\rightarrow0^{+}}{\lim}$
$l_{iL}^{\ast}(\omega)$ exists and $\varsigma\in$ $[-\infty,\ln\frac{V_{0}}%
{D}).$ Suppose $\varsigma$ $\in(-\infty,\ln\frac{V_{0}}{D}).$ Now we apply the
\textquotedblleft continuous dependence\textquotedblright\ theorem of o.d.e.
to the solution of (\ref{2.10a}),(\ref{2.10b}) with the initial value and
parameter pair $((\ln\frac{2V_{0}}{D},0),0).$ When we choose a long enough
existence interval $[0,s_{R}]$ such that the $l$ value at $s_{R}$ of the above
solution is less than $\varsigma-1,$ we have $l_{iL}^{\ast}(\omega
)<\varsigma-1$ as $\omega$ is close to $0$ enough by $\underset{\omega
\rightarrow0^{+}}{\lim}l_{iR}^{\ast}(\omega)=\ln\frac{2V_{0}}{D}$ above and
the \textquotedblleft continuous dependence\textquotedblright\ property. This
contradicts to$\underset{\omega\rightarrow0^{+}}{\lim}$ $l_{iL}^{\ast}%
(\omega)=\varsigma.$ Thus $\underset{\omega\rightarrow0^{+}}{\lim}l_{iL}%
^{\ast}(\omega)=-\infty$ and then (iii) follows.
\end{proof}

\bigskip

Note that by Lemma \ref{Lemma4} and its proof we can roughly say that given
$\omega_{1}<\omega_{2},$ $(\widetilde{l}(s;\omega_{2}),\widetilde{v}%
(s;\omega_{2}))$ traces out an orbit curve which is closer to $(\ln\frac
{V_{0}}{D},0)$ than $(\widetilde{l}(s;\omega_{1}),\widetilde{v}(s;\omega
_{1}))$ does as $s\rightarrow\infty,$ i.e., $(\widetilde{l}(s;\omega
_{2}),\widetilde{v}(s;\omega_{2}))$ is \textquotedblleft
enclosed\textquotedblright\ by $(\widetilde{l}(s;\omega_{1}),\widetilde
{v}(s;\omega_{1})).$

\begin{lemma}
\label{Lemma5}Given any $l(0)\in(-\infty,\ln\frac{2V_{0}}{D}),$ $G,$ and
$i\in\{0,1,2,...\},$ there is at most one $\omega=\omega(G,i,l(0))$ among
$(0,\frac{2V_{0}^{2}}{D})$ such that (\ref{2.10a}),(\ref{2.10b}),(\ref{2.10c})
has a global solution $(l(\cdot;\omega),v(\cdot;\omega))$ on $[0,\infty)$
which has exact $2i$ intersection points with the $l-$axis as $s$ increases
from $0$ to $\infty.$
\end{lemma}

\begin{proof}
Note that any global solution $(l(\cdot;\omega),v(\cdot;\omega))$ of
(\ref{2.10a}),(\ref{2.10b}) has the same orbit curve as $(\widetilde{l}%
(\cdot;\omega),\widetilde{v}(\cdot;\omega)),$ i.e., $(l(s;\omega
),v(s;\omega))=(\widetilde{l}(s+s_{0};\omega),\widetilde{v}(s+s_{0};\omega))$
for some $s_{0}$ by Proposition \ref{Proposition2}. Therefore without loss of
generality, we consider $(\widetilde{l}(\cdot;\omega),\widetilde{v}%
(\cdot;\omega))$ on $[0,\infty)$. Fix $l=$ $l(0)$ and $G$ in (\ref{2.10c}).
Note that the initial value $(l(0),v(0))$ is on the curve:
\begin{equation}
v=I(l;\omega,G):=-\frac{\omega}{D}e^{-l}+\frac{G}{D}. \label{3.10}%
\end{equation}
On one hand $v_{0}:=v(0)=I(l(0);\cdot,G)$ is a decreasing function of $\omega$
on $(0,\infty).$ Let $(\widetilde{l}_{i}(\cdot;\omega),\widetilde{v}_{i}%
(\cdot;\omega))$ be the segment of $(\widetilde{l}(\cdot;\omega),\widetilde
{v}(\cdot;\omega))^{\prime}s$ orbit curve connecting $l_{(i+1)R}^{\ast}%
(\omega)$ and $l_{iL}^{\ast}(\omega)$ in the lower $l-v$ plane. Then on the
other hand, given $\omega_{1}<\omega_{2}$ and $s_{1},s_{2}$ such that
$\widetilde{l}_{i}(s_{1};\omega_{1})=\widetilde{l}_{i}(s_{2};\omega
_{2})=l(0),$ it can be easily seen from Lemma \ref{Lemma4} (i) and its proof
that $\widetilde{v}_{i}(s_{1};\omega_{1})<\widetilde{v}_{i}(s_{2};\omega
_{2}),$ if both $\widetilde{v}_{i}(s_{1};\omega_{1}),\widetilde{v}_{i}%
(s_{2};\omega_{2})$ are nonpositive, i.e., increasing in $\omega.$ Since there
is at most one intersection point for increasing and decreasing functions,
this lemma then follows.
\end{proof}

\bigskip

Now we state our main result:

\begin{theorem}

\begin{enumerate}
\item[(i).] \label{Theorem6}Given any $l(0)\in(-\infty,\ln\frac{2V_{0}}{D})$
and $i\in\{0,1,2,...\},$ there is a corresponding $G_{i}=G_{i}(l(0))>0$ such
that for any given $G\in\lbrack0,G_{i}]$ there is a unique $\omega
=\omega(G;i,l(0))\in(0,\frac{2V_{0}^{2}}{D})$ such that (\ref{2.10a}%
),(\ref{2.10b}),(\ref{2.10c}) has a global existence solution $(l(\cdot
;\omega),v(\cdot;\omega))$ on $(-\infty,\infty)$ which rotates clockwise
around $(\ln\frac{V_{0}}{D},0),$ has $2i$ or $2i+1$ (exact $2i$ for $G=0$
case) intersection points with the $l-$axis as $s$ increases from $0$ to
$\infty,$ and after the $2i$ or $2i+1$ intersections we have $v(s;\omega)<0$,
$l(s;\omega)\searrow-\infty,$ $v(s;\omega)\rightarrow0^{-}$ as
$s\longrightarrow\infty$. Moreover, given any $G\in\underset{i}{\cap}%
[0,G_{i}(l(0))),$ $\omega(G;\cdot,l(0))$ is a decreasing function with the
properties that $\omega(G;i,l(0))\searrow0^{+}$ $as$ $i\rightarrow\infty$ for
$l(0)\in(-\infty,\ln\frac{2V_{0}}{D})\backslash\{\ln\frac{V_{0}}{D}\}$,
$G\in\underset{i}{\cap}[0,G_{i}(l(0)))$ and $\omega(0;i,\ln\frac{V_{0}}%
{D})\searrow0^{+}$ $as$ $i\rightarrow\infty.$

\item[(ii).] For $l(0)\in\lbrack\ln\frac{2V_{0}}{D},\infty)$ and any $G,$ we
have that $\omega=0$ is the only number such that (\ref{2.10a}),(\ref{2.10b}%
),(\ref{2.10c}) has global existence solution $(l(s;0),v(s;0))$. Moreover, if
$G>0,$ then there is some $\vartheta\in(0,\infty)$ such that $v(s;0)>0$ for
$s\in\lbrack0,\vartheta)$, $v(\vartheta;0)=0,$ and $v(s;0)<0$ for
$s\in(\vartheta,\infty)$ with the asymptotic behavior $l(s;0)\searrow-\infty,$
$v(s;0)\rightarrow0^{-}$ as $\vartheta\leq s\longrightarrow\infty;$ if
$G\leq0,$ then $v(s;0)<0$ for $s\in(0,\infty)$ with $l(s;0)\searrow-\infty,$
$v(s;0)\rightarrow0^{-}$ as $0\leq s\longrightarrow\infty.$

\item[(iii).] If $G\in(-\infty,-V_{0}],$ (\ref{2.10a}),(\ref{2.10b}%
),(\ref{2.10c}) has no global existence solutions. When $G\in(-V_{0},0),$
given any $l(0)\in(\ln[\frac{(V_{0}-\sqrt{V_{0}^{2}-G^{2}})}{D}],\ln
[\frac{(V_{0}+\sqrt{V_{0}^{2}-G^{2}})}{D}])$ and $i\in\{0,1,2,...\},$ there is
a unique $\omega=\omega(i,l(0);G)\in(0,\frac{2V_{0}^{2}}{D})$ such that
(\ref{2.10a}),(\ref{2.10b}),(\ref{2.10c}) has a global existence solution
$(l(\cdot;\omega),v(\cdot;\omega))$ satisfying all the properties described in
(i) ( but having exact $2i$ intersection points with the $l-$axis as $s$
increases from $0$ to $\infty$ here); on the other hand, (\ref{2.10a}%
),(\ref{2.10b}),(\ref{2.10c}) has no global existence solutions whenever
$l(0)\notin(\ln[\frac{(V_{0}-\sqrt{V_{0}^{2}-G^{2}})}{D}],\ln[\frac
{(V_{0}+\sqrt{V_{0}^{2}-G^{2}})}{D}])$ in this case$.$
\end{enumerate}
\end{theorem}

First, we explain the physical meaning of the above results and then prove
this theorem.

\subsection{Interpretation of Theorem~\ref{Theorem6}}

(i). Remember that $l(0)=\ln\kappa_{0},$ where $\kappa_{0}$ is the tip's
curvature and $G$ is the tangential velocity of the tip. Theorem
\ref{Theorem6} (i) tells us that given a weakly excitable medium ($V_{0},D$
are fixed), if the tip's curvature $\kappa_{0}$ is not too large, i.e.,
$\kappa_{0}\in(0,\frac{2V_{0}}{D}),$ there is a corresponding range for the
tangential velocity of the tip so that for each $G$ in this range there is a
unique rotating frequency $\omega(G;i,\ln\kappa_{0})$ $\in(0,\frac{2V_{0}^{2}%
}{D})$ such that this medium supports a normally propagating, steadily
rotating, and growing (since $G\geq0$) plane curve with angular frequency
$\omega(G;i,\ln\kappa_{0})$ and tip curvature $\kappa_{0}$, the curvature of
which $\kappa(s)=e^{l(s)}$ changes its monotonicity $2i$ or $2i+1$ times
(exact $2i$ for $G=0$ case) as arc length $s\rightarrow\infty$ and finally
decreases to $0$ in the Archimedean spiral's way, i.e., $\underset
{s\rightarrow\infty}{\lim}\kappa(s)^{2}s=\frac{\omega(G;i,\ln\kappa_{0}%
)}{2V_{0}}$ (see (14) in \cite{IIY})$.$ Especially, the case \textquotedblleft%
$i=0,G=0\textquotedblright$ is just the one appearing in \cite{YaN}.

(ii). On the other hand, if the tip's curvature $\kappa_{0}$ is too large,
i.e., $\kappa_{0}\in\lbrack\frac{2V_{0}}{D},\infty)$, then the medium supports
only \textquotedblleft\emph{nonrotating}\textquotedblright\ spiral waves
($\omega=0$). Moreover, given any $G>0$ ($G\leq0,$ respectively)$,$ there is a
unique growing (contracting, respectively) spiral wave with tip's curvature
$\kappa_{0}$ and tangential velocity $G,$ the curvature of which changes its
monotonicity exactly once ($0$ time, respectively ) and the asymptotic
behavior of which is the same as in above (i) but the curvature $\kappa(s)$
decays to $0$ in different order, not Archimedean one (see THEOREM 2 of
\cite{IIY}).

(iii) tells us that it is impossible for an excitable medium to support a
spiral wave with too large contracting velocity of the tip, i.e.,
$G\in(-\infty,-V_{0}].$ For smaller contracting velocity, $G\in(-V_{0},0),$
there is a corresponding curvature's range of the tip, $(\frac{(V_{0}%
-\sqrt{V_{0}^{2}-G^{2}})}{D},\frac{(V_{0}+\sqrt{V_{0}^{2}-G^{2}})}{D}),$ such
that when the tip's curvature $\kappa_{0}$ does not belong to this range,
there is no any spiral wave but when given any $\kappa_{0}$ in this range,
there is a unique $\omega(i,\ln\kappa_{0};G)\in(0,\frac{2V_{0}^{2}}{D})$ such
that this medium supports a spiral wave with tip's contracting velocity $G,$
tip's curvature $\kappa_{0},$ rotating frequency $\omega(i,\ln\kappa_{0};G),$
and the same properties as mentioned in (i) above (but the curvature changes
its monotonicity exact $2i$ times).

Note that we also have similar explanations as above by replacing $\kappa,$
$\omega,$ $V_{0}$ with $-\kappa,$ $-\omega,$ $-V_{0},$ respectively when
$\kappa<0$.

\subsection{Proof of Theorem~\ref{Theorem6}}

\begin{proof}
As explained in the proof of Lemma \ref{Lemma5}, it suffices to consider
$(\widetilde{l}(\cdot;\omega),\widetilde{v}(\cdot;\omega)).$ Letting
$\omega=\frac{2V_{0}^{2}}{D},$ (\ref{3.10}) becomes
\begin{equation}
v=I(l;\frac{2V_{0}^{2}}{D},G):=-2(\frac{V_{0}}{D})^{2}e^{-l}+\frac{G}{D}.
\label{3.11}%
\end{equation}
Plugging (\ref{3.11}) into $E(l,v),$ we obtain
\begin{equation}
E(l,-2(\frac{V_{0}}{D})^{2}e^{-l}+\frac{G}{D})=\frac{1}{2}e^{-2l}%
[4(\frac{V_{0}}{D})^{4}+e^{4l}-\frac{2V_{0}}{D}e^{3l}+(\frac{G}{D})^{2}%
e^{2l}-4(\frac{G}{D})(\frac{V_{0}}{D})^{2}e^{l}]. \label{3.12}%
\end{equation}
It can be easily derived that \textquotedblleft$x^{4}-2ax^{3}+4a^{4}>0,$
$\forall x\in R,a\neq0$\textquotedblright\ which implies $E(l,-2(\frac{V_{0}%
}{D})^{2}e^{-l})>0$ for all $l\in R$ (i.e., when $G=0$)$.$ Hence given any
$l,$ we have that $I(l;\frac{2V_{0}^{2}}{D},0)<$the negative $v$ value of
curve \textquotedblleft$E(l,v)=0\textquotedblright,$ i.e., curve (\ref{3.11})
lies below $E(l,v)=0$ when $G=0.$ Then by continuity, there is a $G(l)>0$ such
that $I(l;\frac{2V_{0}^{2}}{D},G)$ $\leq$ the negative $v$ value of curve
\textquotedblleft$E(l,v)=0\textquotedblright,$ whenever $G\in(-\infty,G(l)].$
From above observation and Lemma \ref{Lemma0.2}, we have that if there is some
$\omega>0$ such that (\ref{2.10a}),(\ref{2.10b}),(\ref{2.10c}) has global
solutions, then $\omega\in(0,\frac{2V_{0}^{2}}{D}).$ Thus we have obtained the
uniqueness property by Lemma \ref{Lemma5}. It is left for (i) to prove the
\textquotedblleft existence\textquotedblright\ of such $\omega.$ The strategy
is to find suitable $\omega$ such that there are intersections points between
the curve (\ref{3.10}) and some special segment of the orbit curve of
$(\widetilde{l}(\cdot;\omega),\widetilde{v}(\cdot;\omega)),$ which implies
(\ref{2.10a}),(\ref{2.10b}),(\ref{2.10c}) has global solutions with the
properties mentioned in this theorem$.$ By Lemma \ref{Lemma4}, the ranges of
$l_{(i+1)R}^{\ast}(\cdot)$ and $l_{iL}^{\ast}(\cdot)$ over $(0,\frac
{2V_{0}^{2}}{D})$ are open intervals, say, $(R_{i+1},\ln\frac{2V_{0}}{D})$ and
$(-\infty,L_{i})$ respectively, where $l_{0L}^{\ast}(\cdot):=-\infty$ and
$L_{0}:=-\infty.$ For any $i\in\{0,1,2,...\}$ and $l(0)\in(L_{i},R_{i+1}),$
$\widetilde{v}_{i}(l(0);\omega)$ is well-defined for any $\omega\in
(0,\frac{2V_{0}^{2}}{D}),$ where $\widetilde{v}_{i}(l;\omega)$ $:=\widetilde
{v}(s;\omega)$ if there is some $s$ such that $\widetilde{l}(s;\omega)=l$ and
$(\widetilde{l}(s;\omega),\widetilde{v}(s;\omega))$ is on the
\textquotedblleft lower\textquotedblright\ orbit curve connecting
$l_{(i+1)R}^{\ast}(\omega)$ and $l_{iL}^{\ast}(\omega).$ Thus given any
$l(0)\in(L_{i},R_{i+1})$ and $G\in(-\infty,G_{i}(l(0))],$ where $G_{i}%
(\cdot):=G(\cdot),$ there is some positive $\omega_{P_{i}}$ close to
$\frac{2V_{0}^{2}}{D}$ such that $\widetilde{v}_{i}(l(0);\omega_{P_{i}})>$the
negative $v$ value of \textquotedblleft$E(l(0),v)=0\textquotedblright\geq
I(l(0);\omega_{P_{i}},G)$ by Lemma \ref{Lemma0.2}, the continuity of
$I(l(0);\cdot,G),$ and the observation in the beginning of this proof; for
$l(0)=L_{i}$ or $R_{i+1},$ it still holds by Lemma \ref{Lemma4} (i),(ii). Now
consider $l(0)\in(R_{i+1},\ln\frac{2V_{0}}{D})$ (or $l(0)\in(-\infty,L_{i})).$
By definition, there is $\omega_{0}\in(0,\frac{2V_{0}^{2}}{D})$ such that
$l_{(i+1)R}^{\ast}(\omega_{0})=l(0)$ (or $l_{iL}^{\ast}(\omega_{0})=l(0)).$ We
require now that $I(l(0);\omega_{0},G)\leq0$ which is equivalent to
\begin{equation}
G\leq\omega_{0}e^{-l(0)}. \label{3.13}%
\end{equation}
Let $G_{i}(l):=\omega_{0}e^{-l}.$ Therefore, given any $G\in(-\infty
,G_{i}(l(0))],$ there is also some $\omega_{P_{i}}$ close to $\omega_{0}$ such
that $\widetilde{v}_{i}(l(0);\omega_{P_{i}})\geq$ $I(l(0);\omega_{P_{i}},G)$
by virture of Lemma \ref{Lemma4} (i),(ii). On the other hand, given any
$G\in\lbrack0,\infty)$ and $l(0)\in(-\infty,\ln\frac{2V_{0}}{D}),$ there is
some positive $\omega_{N_{i}}$ close to $0$ such that $\widetilde{v}%
_{i}(l(0);\omega_{N_{i}})$ is well-defined and $\widetilde{v}_{i}%
(l(0);\omega_{N_{i}})<$ $I(l(0);\omega_{N_{i}},G),$ since $\underset
{\omega\rightarrow0^{+}}{\lim}$ $I(l(0);\omega,G)=\frac{G}{D}\geq0$ and
$\underset{\omega\rightarrow0^{+}}{\lim}\widetilde{v}_{i}(l(0);\omega)=$the
negative $v$ value of \ \textquotedblleft$E(l(0),v)=0"<0$ by Lemma
\ref{Lemma4} (iii) and the continuous dependence theorem of o.d.e. Note that
$\widetilde{v}_{i}(l(0);\cdot)$ is a continuous function on $[\omega_{N_{i}%
},\omega_{P_{i}}]$ by Lemma \ref{Lemma4} (ii)$.$ Obviously $I(l(0);\cdot,G)$
is also a continuous function.\ Therefore choosing $l(0)\in(-\infty,\ln
\frac{2V_{0}}{D})$ and then $G\in\lbrack0,G_{i}(l(0))),$ there must be some
$\sigma_{i}\in(\omega_{N_{i}},\omega_{P_{i}}]\subset(0,\frac{2V_{0}^{2}}{D})$
such that $\widetilde{v}_{i}(l(0);\sigma_{i})=I(l(0);\sigma_{i},G)$ by the
intermediate value theorem of continuous functions. Let $\omega
(G;i,l(0)):=\sigma_{i}.$ Hence \textquotedblleft existence\textquotedblright%
\ is guaranteed. Given $G\in\underset{i}{\cap}[0,G_{i}(l(0))),$ the
decreasingness of $\omega(G;\cdot,l(0))$ as $i\rightarrow\infty$ immediately
follows from the explanation after Lemma \ref{Lemma4}, i.e., \textquotedblleft%
$(\widetilde{l}(s;\omega_{2}),\widetilde{v}(s;\omega_{2}))$ is enclosed\ by
$(\widetilde{l}(s;\omega_{1}),\widetilde{v}(s;\omega_{1}))$ as $\omega
_{1}<\omega_{2}$\textquotedblright\ and the decreasingness of $I(l(0);\cdot
,0)$. Now we prove $\underset{i\rightarrow\infty}{\lim}\omega(G;i,l(0))=0$ for
$l(0)\in(-\infty,\ln\frac{2V_{0}}{D})\backslash\{\ln\frac{V_{0}}{D}\}$,
$G\in\underset{i}{\cap}[0,G_{i}(l(0)))$ and $\underset{i\rightarrow\infty
}{\lim}\omega(0;i,\ln\frac{V_{0}}{D})=0.$ For the first part, assume the
contrary that $\underset{i\rightarrow\infty}{\lim}\omega(G;i,l(0))=\omega
_{f}(G;l(0))>0.$ Note that $\underset{i\rightarrow\infty}{\lim}l_{iR}^{\ast
}(\omega_{f}(G;l(0)))=\underset{i\rightarrow\infty}{\lim}l_{iL}^{\ast}%
(\omega_{f}(G;l(0)))=\ln\frac{V_{0}}{D}.$ By $l(0)\neq\ln\frac{V_{0}}{D},$
there must be some $i$ large enough such that $l(0)\in(l_{(i+1)R}^{\ast
}(\omega_{f}(G;l(0))),\ln\frac{2V_{0}}{D})$ $\subset(R_{i+1},\ln\frac{2V_{0}%
}{D})$ or $l(0)\in(-\infty,l_{iL}^{\ast}(\omega_{f}(G;l(0))))\subset
(-\infty,L_{i}).$ Arguing as before, then there is some $\omega\in
(0,\omega_{f}(G;l(0)))$ such that $\widetilde{v}_{i}(l(0);\omega
)=I(l(0);\omega,G).$ By uniqueness, this $\omega$ must be $\omega(G;i,l(0))$
which leads to a contradiction, since it is impossible that $\omega
(G;i,l(0))<$ $\omega_{f}(G;l(0)).$ For the second part, note that
$\widetilde{v}_{i}(\ln\frac{V_{0}}{D};\omega_{f}(0;\ln\frac{V_{0}}{D}))$ is
always defined and $\underset{i\rightarrow\infty}{\lim}\widetilde{v}_{i}%
(\ln\frac{V_{0}}{D};\omega_{f}(0;\ln\frac{V_{0}}{D}))=0$ as explained in the
last of the proof of Lemma \ref{Lemma1}. Thus for $i$ large enough we have
$\widetilde{v}_{i}(\ln\frac{V_{0}}{D};\omega_{f}(0;\ln\frac{V_{0}}{D}%
))>I(\ln\frac{V_{0}}{D};\omega_{f}(0;\ln\frac{V_{0}}{D}),0),$ since
$I(\ln\frac{V_{0}}{D};\omega_{f}(0;\ln\frac{V_{0}}{D}),0)<0.$ Then again by
Lemma \ref{Lemma4} (iii) and the intermediate value theorem, there must be
some $\omega\in(0,\omega_{f}(0;\ln\frac{V_{0}}{D}))$ such that $\widetilde
{v}_{i}(\ln\frac{V_{0}}{D};\omega)=I(\ln\frac{V_{0}}{D};\omega,0)$ which also
leads to a contradiction$.$This completes (i).

Part (ii) follows from Lemma \ref{Lemma0.2} and THEOREM 2 in \cite{IIY} immediately.

Now we prove (iii). When $G\in(-\infty,-V_{0}],$ the curve $(l,I(l;\omega,G))$
entirely lies outside the curve \textquotedblleft$E(l,v)=0$\textquotedblright%
\ for any $\omega>0$. Thus by Lemma \ref{Lemma0.2}, we have that
(\ref{2.10a}),(\ref{2.10b}),(\ref{2.10c}) has no global existence solutions.
When $G\in(-V_{0},0),$ the parts of the curve $(l,I(l;\omega,G)),$ where
$l(0)\notin(\ln[\frac{(V_{0}-\sqrt{V_{0}^{2}-G^{2}})}{D}],\ln[\frac
{(V_{0}+\sqrt{V_{0}^{2}-G^{2}})}{D}]),$ also lie outside the curve
\textquotedblleft$E(l,v)=0$\textquotedblright\ for any $\omega>0.$ As
explained above, we only need to consider $l(0)\in(\ln[\frac{(V_{0}%
-\sqrt{V_{0}^{2}-G^{2}})}{D}],\ln[\frac{(V_{0}+\sqrt{V_{0}^{2}-G^{2}})}{D}])$
now. Note that the right-hand side of (\ref{3.12}) is greater than $0$
$\forall$ $l,$ for any $G<0$ and (\ref{3.13}) is also valid. Arguing as in
(i), we then have completed (iii).

Therefore we have finished the proof.
\end{proof}

\bigskip

\section{Numerical Results}

\label{numres}

The following formulas can be referred to \cite[p. 45]{Mer}.

Let
\begin{equation}
\theta_{0}(t)=\omega t+\theta_{0}(0) \label{4.1}%
\end{equation}
be the angle of the tip at time $t$, where $\omega$ is the rotating frequency
determined by Theorem~\ref{Theorem6} and $\theta_{0}(0)$ is the initial angle
of the tip, which can be chosen arbitrary. Let also $\theta(s,t)$ be the angle
of the position away from the tip at a distance \textquotedblleft arclength
$s$\textquotedblright\ at time $t.$ Then
\begin{equation}
\theta(s,t)=\theta_{0}(t)-\int_{0}^{s}e^{l(\rho)}d\rho, \label{4.2}%
\end{equation}
where $l(s)$ is the function from Theorem~\ref{Theorem6}.

The cartesian coordinates of the tip at time $t$ are given by
\begin{align}
\frac{d}{dt}X_{0}(t)  &  =(De^{l(0)}-V_{0})\sin\theta_{0}(t)-G\cos\theta
_{0}(t)=(De^{l(0)}-V_{0})\sin[\omega t+\theta_{0}(0)]-G\cos[\omega
t+\theta_{0}(0)]\quad\text{and}\label{4.3}\\
\frac{d}{dt}Y_{0}(t)  &  =\left(  V_{0}-De^{l(0)}\right)  \cos\theta
_{0}(t)-G\sin\theta_{0}(t)=\left(  V_{0}-De^{l(0)}\right)  \cos[\omega
t+\theta_{0}(0)]-G\sin[\omega t+\theta_{0}(0)],\nonumber
\end{align}
where the initial condition $(X_{0}(0),\,Y_{0}(0))$ can be chosen arbitrarily.
Then, the cartesian coordinates $(X(s,t),Y(s,t))$ of the position away from
the tip at a distance \textquotedblleft arclength $s$\textquotedblright\ at
time $t$ are
\begin{align}
X(s,t)  &  =\int_{0}^{s}\cos\theta(\xi,t)d\xi+X_{0}(t)=\int_{0}^{s}\cos[\omega
t+\theta_{0}(0)-\int_{0}^{\xi}e^{l(\rho)}d\rho]d\xi+X_{0}(t)\label{4.4}\\
Y(s,t)  &  =\int_{0}^{s}\sin\theta(\xi,t)d\xi+Y_{0}(t)=\int_{0}^{s}\sin[\omega
t+\theta_{0}(0)-\int_{0}^{\xi}e^{l(\rho)}d\rho]d\xi+Y_{0}(t),\nonumber
\end{align}
which is the same as (4.43) in \cite[p. 45]{Mer}.

Now we explain how to determine the parameters $D,$ $V_{0},$ $\theta_{0}(0),$
$l(0),$ $G,$ $\omega$ and then to draw the pictures of spiral waves \ in
Theorem~\ref{Theorem6}.

First, parameters $D>0,$ $V_{0}>0,$ $\theta_{0}(0)$ can be chosen arbitrarily.

In Theorem~\ref{Theorem6} (i), we choose $l(0)\in(-\infty,\ln\frac{2V_{0}}%
{D}),$ $i\in\{0,1,2,3,...\}$ first$.$ Then we choose $G\in\lbrack
0,G_{i}(l(0))]$ (for example choose $G=0$). Finally, there is a unique
$\omega=\omega(G;i,l(0))\in(0,\frac{2V_{0}^{2}}{D})$ by Theorem~\ref{Theorem6}
(i). Plug above parameters into (\ref{4.3}) and (\ref{4.4})\ and then the
corresponding spiral wave has the properties described in
Theorem~\ref{Theorem6} (i).

For (ii), we choose $l(0)\in\lbrack\ln\frac{2V_{0}}{D},\infty),$ $\omega=0,$
and any $G.$ Plug above parameters into (\ref{4.3}) and (\ref{4.4})\ and then
the corresponding spiral wave has the properties described in
Theorem~\ref{Theorem6} (ii).

For (iii), we first choose $G\in(-V_{0},0)$ and then $l(0)\in(\ln[\frac
{(V_{0}-\sqrt{V_{0}^{2}-G^{2}})}{D}],\ln[\frac{(V_{0}+\sqrt{V_{0}^{2}-G^{2}}%
)}{D}])$ and any $i\in\{0,1,2,3,...\}.$ Finally there is a unique
$\omega(i,l(0);G)\in(0,\frac{2V_{0}^{2}}{D}).$ Also, plugging above parameters
into (\ref{4.3}) and (\ref{4.4}),\ the corresponding spiral wave then has the
properties described in Theorem~\ref{Theorem6} (iii).

\medskip

\paragraph{\textbf{Acknowledgements}}

C.-P.Lo thanks Mr. J.-C.Tsai and professors J.-S.Guo and S.-S.Lin for helpful suggestions.

\end{document}